\documentclass[12pt]{amsart}
\usepackage[all]{xy}

\newtheorem{thm}{Theorem}[section]
\newtheorem{prob}[thm]{Problem}

\newcommand{\impl}{\rightarrow}
\newcommand{\w}{\omega}
\renewcommand{\b}{\mathfrak{b}}

\renewcommand{\split}{\mathsf{Split}}
\newcommand{\bq}{\begin{quote}}
\newcommand{\eq}{\end{quote}}
\renewcommand{\O}{\mathcal{O}}
\newcommand{\B}{\mathcal{B}}
\newcommand{\BG}{\B_\Gamma}

\newcommand{\BO}{\B_\Omega}

\newcommand{\sone}{\mathsf{S}_1}    \newcommand{\sfin}{\mathsf{S}_{fin}}
\newcommand{\ufin}{\mathsf{U}_{fin}}
    
\newcommand{\seq}[1]{\{#1\}_{n\in\N}}
\newcommand{\Union}{\bigcup}

\newcommand{\fU}{\mathfrak{U}}
\newcommand{\fV}{\mathfrak{V}}

\newcommand{\naturals}{{\mathbb N}}
\newcommand{\N}{\naturals}

\newcommand{\by}[2]{\hfill\emph{#1}, #2}

\newcommand{\CE}{\textsc{ce}}

\newcommand{\be}{\begin{enumerate}}
\newcommand{\ee}{\end{enumerate}}
\newcommand{\bi}{\begin{itemize}}
\newcommand{\ei}{\end{itemize}}
\renewcommand{\i}{\item}
\newcommand{\SPMBul}{\textbf{$\mathcal{SPM}$ BULLETIN}}

\newcommand{\arx}[1]{\texttt{http://arxiv.org/abs/#1}}

\title[\SPMBul{} \textbf{2} (February 2003)]{%
\SPMBul\\[0.5cm]
Issue number 2: February 2003 \CE{}}

\begin{document}
\maketitle

\tableofcontents

\section{Editor's note}

This issue of the \SPMBul{} contains
several new research announcements.
We also have a conference announcement.
This conference seems particularly interesting
for recipients of this bulletin:
Liljana Babinkostova and Gary Gruenhage are
experts in the field of topological diagonalizations,
Klaas Pieter Hart is an expert in filters, and
Jindrich Zapletal is an expert in cardinal characteristics of the continuum.

Recall that in most cases, the information does not reach us
without your help. Please contact us with any interesting
information you have on new preprints/papers,
conferences/workshops, and candidates for subscription
(our list of recipients expanded a little in the last month,
but we know that there are many more mathematicians working in the field
which may be interested in receiving this bulletin).

The first issue of this bulletin is available online at
\bq
\arx{math.GN/0301011}
\eq
In that issue we gave some of the basic definitions used frequently in the field
of Selection Principles in Mathematics.
In the \emph{Problem of the month} section we describe
the classification of the resulting properties with regards
to (standard) covers, large covers, $\w$-covers, and $\gamma$-covers.
Only two classification problems remain unsettled and constitute the
problem of the month.

\subsection{Dedication}
After this issue was completed, we heard of the
tragedy of space shuttle Columbia, who
broke up over Texas on Feb 1 2003 \CE{}, as it headed for a landing in Florida.
We dedicate this issue to the memorial
of the shuttle crew:
Michael Anderson, David Brown, Kalpana Chawla, Laurel Clark,
Rick Husband, William McCool, and Ilan Ramon.

\by{Boaz Tsaban}{tsaban@math.huji.ac.il}

\hfill \texttt{http://www.cs.biu.ac.il/\~{}tsaban}

\section{Research announcements}

\subsection{The Pytkeev property and the Reznichenko property in function spaces}
For a Tychonoff space $X$
we denote by $C_p(X)$
the space of all real-valued continuous functions on $X$ with
the topology of pointwise convergence.
Characterizations of sequentiality and countable tightness
of $C_p(X)$ in terms of $X$ were given by
Gerlits, Nagy, Pytkeev and Arhangel'skii.
In this paper, we characterize
the Pytkeev property and the Reznichenko property
of $C_p(X)$ in terms of $X$.
In particular we note that
if $C_p(X)$ over a subset $X$ of the real line is a Pytkeev space,
then $X$ is perfectly meager and has universal measure zero.

\by{Masami Sakai}{sakaim01@kanagawa-u.ac.jp}

\subsection{The combinatorics of splittability}
Marion Scheepers, in his studies of the combinatorics of open covers,
introduced the property $\split(\fU,\fV)$ asserting that a cover of
type $\fU$ can be split into two covers of type $\fV$.
In the first part of this paper
we give an almost complete classification of all properties of this
form where $\fU$ and $\fV$
are significant families of covers which appear in the literature
(namely, large covers, $\omega$-covers, $\tau$-covers, and $\gamma$-covers),
using combinatorial characterizations of these properties in terms related
to ultrafilters on $\N$.

In the second part of the paper we
consider the questions whether, given $\fU$ and $\fV$,
the property $\split(\fU,\fV)$ is preserved under taking finite unions,
arbitrary subsets, powers or products.
Several interesting problems remain open.

The paper is available online:
\bq
\arx{math.LO/0212312}
\eq

\by{Boaz Tsaban}{tsaban@math.huji.ac.il}

\subsection{The Hurewicz covering property and slaloms in the Baire space}

According to a result of Ko\v{c}inac and Scheepers,
the Hurewicz covering property is equivalent to a somewhat simpler
selection property: For each sequence of large open covers
of the space one can choose finitely many elements from each cover
to obtain a groupable cover of the space.
We simplify the characterization further by omitting the need to
consider sequences of covers:
A set of reals $X$ satisfies the Hurewicz property if, and only if,
each large open cover of $X$ contains a groupable subcover.

The proof uses a ``structure'' counterpart of a combinatorial
characterization, in terms of slaloms,
of the minimal cardinality $\b$ of an unbounded family
of functions in the Baire space.
In particular, we obtain a new characterization of $\b$.

The paper is available online:
\bq
\arx{math.GN/0301085}
\eq

\by{Boaz Tsaban}{tsaban@math.huji.ac.il}

\subsection{Games in Logic}
A nice review of Games in Logic (including Game-Theoretic
Semantics), written in 2001 by the logician Wilfrid Hodges for the
Stanford Encyclopedia of Philosophy, is available online at:
\verb|http://plato.stanford.edu/entries/logic-games/|

\by{Peter McBurney}{p.j.mcburney@csc.liv.ac.uk}

\section{Other announcements}

\subsection{Boise Extravaganza in Set Theory (BEST 2003)}
Friday, March 28 - Sunday, March 30, 2003.

We are pleased to announce our twelfth annual BEST conference.
There will be 4 talks by invited speakers:
\bi
\i Liljana Babinkostova (Boise State University)
\i Klaas Pieter Hart (Delft University of Technology, the Netherlands)
\i Gary Gruenhage (Auburn University)
\i Jindrich Zapletal (University of Florida)
\ei

The talks will be held on Friday, Saturday and on Sunday in the
Department of Mathematics at Boise State University.
BEST social events are planned as well.

The conference webpage at
\bq
            \verb|http://math.boisestate.edu/~best/|
\eq
contains the most current information including lodging, abstract
submission, maps, schedule, etc.
Anyone interested in giving a talk and/or participating
should contact one of the organizers as soon as possible.

Limited financial support is available, in particular for graduate
students. In order to apply, e-mail one of the organizers and we will
send you the required paperwork.

The conference is  supported by a grant from the National
Science Foundation, whose assistance is gratefully acknowledged.

\by{Tomek Bartoszynski}{tomek@math.boisestate.edu}

\by{Justin Moore}{justin@math.boisestate.edu}

\subsection{A web site dedicated to SPM}
Marion Scheepers and Liljana Babinkostova are building a comprehensive
web site dedicated to Selection Principles in Mathematics.
According to Scheepers, this site is planned to
``grow to contain everything under the sun related to SPM''.
Currently, the site contains:
\bi
\i A useful bibliography of works in the field,
\i Online versions of some of the major papers in the field,
\i Seminar meetings and workshop announcements; and
\i Selected open problems (under construction).
\ei
The site's address is:
\bq
\verb|http://iunona.pmf.ukim.edu.mk/~spm/|
\eq

\by{Boaz Tsaban}{tsaban@math.huji.ac.il}

\section{Problem of the month}
The following discussion is borrowed from \cite{coc2, CBC}.
(The paper \cite{coc2} makes an excellent introduction to
one of the major aspects of the field of SPM.)

We will use the notation given in the first issue of this
bulletin (\arx{math.GN/0301011}).
For the types of covers which we consider,
$$\sone(\fU,\fV)\impl\sfin(\fU,\fV)\impl\ufin(\fU,\fV)\&\binom{\fU}{\fV}$$
and $\binom{\Lambda}{\Omega}$ does not hold for an infinite $T_1$ space $X$
\cite{coc2, strongdiags}. This rules out several of the introduced properties
as trivial.
Each of our properties is monotone decreasing in the first coordinate
and increasing in the second.
In the case of $\ufin$ note that
for any class of covers $\fV$,
$\ufin(\O,\fV)$ is equivalent to $\ufin(\Gamma,\fV)$
because given an open cover $\seq{U_n}$ we may replace
it by $\seq{\Union_{i<n} U_i}$, which is a $\gamma$-cover
(unless it contains a finite subcover).

In the three-dimensional diagram of Figure~\ref{3dim} below,
the double lines indicate that the two properties are equivalent.
The proof of these equivalences can be found in \cite{coc1, coc2}.

\begin{figure}[!hb]
\unitlength=.8mm
\begin{picture}(141.00,145.00)(0,0)
\put(104.00,20.00){\makebox(0,0)[cc]{$\sone(\O,\O)$}}
\put(104.00,46.00){\makebox(0,0)[cc]{$\sone(\Lambda,\O)$}}
\put(73.00,46.00){\makebox(0,0)[cc]{$\sone(\Lambda,\Lambda)$}}
\put(104.00,72.00){\makebox(0,0)[cc]{$\sone(\Omega,\O)$}}
\put(73.00,72.00){\makebox(0,0)[cc]{$\sone(\Omega,\Lambda)$}}
\put(43.00,72.00){\makebox(0,0)[cc]{$\sone(\Omega,\Omega)$}}
\put(13.00,72.00){\makebox(0,0)[cc]{$\sone(\Omega,\Gamma)$}}
\put(104.00,99.00){\makebox(0,0)[cc]{$\sone(\Gamma,\O)$}}
\put(73.00,99.00){\makebox(0,0)[cc]{$\sone(\Gamma,\Lambda)$}}
\put(43.00,99.00){\makebox(0,0)[cc]{$\sone(\Gamma,\Omega)$}}
\put(13.00,99.00){\makebox(0,0)[cc]{$\sone(\Gamma,\Gamma)$}}

\put(123.00,33.00){\makebox(0,0)[cc]{$\sfin(\O,\O)$}}
\put(123.00,59.00){\makebox(0,0)[cc]{$\sfin(\Lambda,\O)$}}
\put(89.00,59.00){\makebox(0,0)[cc]{$\sfin(\Lambda,\Lambda)$}}
\put(123.00,86.00){\makebox(0,0)[cc]{$\sfin(\Omega,\O)$}}
\put(89.00,86.00){\makebox(0,0)[cc]{$\sfin(\Omega,\Lambda)$}}
\put(58.00,86.00){\makebox(0,0)[cc]{$\sfin(\Omega,\Omega)$}}
\put(27.00,86.00){\makebox(0,0)[cc]{$\sfin(\Omega,\Gamma)$}}
\put(27.00,115.00){\makebox(0,0)[cc]{$\sfin(\Gamma,\Gamma)$}}
\put(58.00,115.00){\makebox(0,0)[cc]{$\sfin(\Gamma,\Omega)$}}
\put(89.00,115.00){\makebox(0,0)[cc]{$\sfin(\Gamma,\Lambda)$}}
\put(123.00,115.00){\makebox(0,0)[cc]{$\sfin(\Gamma,\O)$}}

\put(141.00,130.00){\makebox(0,0)[cc]{$\ufin(\Gamma,\O)$}}
\put(43.00,130.00){\makebox(0,0)[cc]{$\ufin(\Gamma,\Gamma)$}}
\put(73.00,130.00){\makebox(0,0)[cc]{$\ufin(\Gamma,\Omega)$}}
\put(104.00,130.00){\makebox(0,0)[cc]{$\ufin(\Gamma,\Lambda)$}}

\put(105.00,13.00){\makebox(0,0)[cc]{$C^{\prime\prime}$
Rothberger}}
\put(11.00,64.00){\makebox(0,0)[cc]{$\gamma$-set Gerlits-Nagy}}
\put(131.00,137.00){\makebox(0,0)[cc]{Menger}}
\put(43.00,137.00){\makebox(0,0)[cc]{Hurewicz}}

\put(054.00,130.00){\vector(1,0){9.00}}
\put(084.00,130.00){\vector(1,0){9.00}}
\put(114.00,130.50){\line(1,0){16.00}}
\put(114.00,129.50){\line(1,0){16.00}}

\put(061.00,118.00){\vector(1,1){10.00}}
\put(092.00,118.00){\line(1,1){10.00}}
\put(091.00,118.00){\line(1,1){10.00}}
\put(030.00,118.00){\vector(1,1){10.00}}
\put(122.00,118.00){\line(1,1){10.00}}
\put(123.00,118.00){\line(1,1){10.00}}

\put(036.00,115.00){\vector(1,0){12.00}}
\put(067.00,115.00){\vector(1,0){12.00}}
\put(099.00,115.50){\line(1,0){14.00}}
\put(099.00,114.50){\line(1,0){14.00}}

\put(014.50,102.00){\line(1,1){10.00}}
\put(015.50,102.00){\line(1,1){10.00}}
\put(043.00,102.00){\vector(1,1){10.00}}
\put(076.00,102.00){\vector(1,1){10.00}}
\put(107.00,102.00){\vector(1,1){10.00}}

\put(090.00,101.00){\line(0,1){10.00}}
\put(089.00,101.00){\line(0,1){10.00}}
\put(058.00,101.00){\vector(0,1){10.00}}
\put(027.00,101.00){\vector(0,1){10.00}}

\put(021.00,099.00){\vector(1,0){14.00}}
\put(051.00,099.00){\vector(1,0){14.00}}
\put(081.00,099.50){\line(1,0){14.00}}
\put(081.00,098.50){\line(1,0){14.00}}

\put(027.00,090.00){\line(0,1){6.00}}
\put(058.00,090.00){\line(0,1){6.00}}
\put(090.00,090.00){\line(0,1){6.00}}
\put(089.00,090.00){\line(0,1){6.00}}
\put(123.00,090.00){\line(0,1){22.00}}
\put(122.00,090.00){\line(0,1){22.00}}

\put(037.00,086.00){\line(1,0){4.00}}
\put(045.00,086.00){\vector(1,0){3.00}}
\put(068.00,086.00){\line(1,0){3.00}}
\put(074.00,086.00){\vector(1,0){5.00}}
\put(099.00,086.50){\line(1,0){3.00}}
\put(099.00,085.50){\line(1,0){3.00}}
\put(106.00,086.50){\line(1,0){7.00}}
\put(106.00,085.50){\line(1,0){7.00}}

\put(016.00,076.00){\line(1,1){7.00}}
\put(017.00,076.00){\line(1,1){7.00}}
\put(047.00,076.00){\vector(1,1){7.00}}
\put(077.00,076.00){\vector(1,1){7.00}}
\put(107.00,076.00){\vector(1,1){7.00}}

\put(013.00,076.00){\vector(0,1){19.00}}
\put(043.00,076.00){\vector(0,1){19.00}}
\put(072.00,076.00){\vector(0,1){20.00}}
\put(090.00,075.00){\line(0,1){7.00}}
\put(089.00,075.00){\line(0,1){7.00}}
\put(104.00,076.00){\vector(0,1){20.00}}

\put(021.00,072.00){\vector(1,0){14.00}}
\put(051.00,072.00){\vector(1,0){14.00}}
\put(081.00,072.50){\line(1,0){14.00}}
\put(081.00,071.50){\line(1,0){14.00}}

\put(089.00,062.00){\line(0,1){7.00}}
\put(090.00,062.00){\line(0,1){7.00}}
\put(122.00,064.00){\line(0,1){19.00}}
\put(123.00,064.00){\line(0,1){19.00}}

\put(099.00,059.50){\line(1,0){3.00}}
\put(099.00,058.50){\line(1,0){3.00}}
\put(106.00,059.50){\line(1,0){6.00}}
\put(106.00,058.50){\line(1,0){6.00}}

\put(078.00,050.00){\vector(1,1){7.00}}
\put(108.00,050.00){\vector(1,1){7.00}}

\put(073.00,050.00){\line(0,1){19.00}}
\put(072.00,050.00){\line(0,1){19.00}}
\put(104.00,049.00){\line(0,1){19.00}}
\put(103.00,049.00){\line(0,1){19.00}}

\put(082.00,046.00){\line(1,0){13.00}}
\put(082.00,045.00){\line(1,0){13.00}}

\put(123.00,035.00){\line(0,1){19.00}}
\put(122.00,035.00){\line(0,1){19.00}}

\put(104.00,023.00){\line(0,1){19.00}}
\put(103.00,023.00){\line(0,1){19.00}}

\put(108.00,023.00){\vector(1,1){8.00}}

\end{picture}
\caption{\label{3dim}}
\end{figure}
The analogue equivalences for the Borel case also hold, but
in the Borel case more equivalences hold \cite{CBC}:
For each $\fV \in\{\B,\BO,\BG\}$,
$$\sone(\BG,\fV) = \sfin(\BG,\fV) = \ufin(\BG,\fV).$$
After removing duplications we obtain Figure \ref{survive} (see below).

\begin{figure}
\newcommand{\sr}[2]{#1}
{\tiny
$$\xymatrix@C=-2pt@R=10pt{
&
&
& \sr{\ufin(\Gamma,\Gamma)}{\b}\ar[rr]\ar@{.>}[dr]^?
&
& \sr{\ufin(\Gamma,\Omega)}{\d}\ar[rrrrr]\ar@/_/@{.>}[dl]_?
&
&
&
&
&
&
& \sr{\ufin(\Gamma,\O)}{\d}
\\
&
&
&
& \sr{\sfin(\Gamma,\Omega)}{\d}\ar[ur]
\\
& \sr{\sone(\Gamma,\Gamma)}{\b}\ar[rr]\ar[uurr]
&
& \sr{\sone(\Gamma,\Omega)}{\d}\ar[rrr]\ar[ur]
&
&
& \sr{\sone(\Gamma,\O)}{\d}\ar[uurrrrrr]
\\
  \sr{\sone(\BG,\BG)}{\b}\ar[ur]\ar[rr]
&
& \sr{\sone(\BG,\BO)}{\d}\ar[ur]\ar[rrr]
&
&
& \sr{\sone(\BG,\B)}{\d}\ar[ur]
\\
&
&
&
& \sr{\sfin(\Omega,\Omega)}{\d}\ar'[u]'[uu][uuu]
\\
\\
&
& \sr{\sfin(\BO,\BO)}{\d}\ar[uuu]\ar[uurr]
\\
& \sr{\sone(\Omega,\Gamma)}{\p}\ar'[r][rr]\ar'[uuuu][uuuuu]
&
& \sr{\sone(\Omega,\Omega)}{\cov(\M)}\ar'[uuuu][uuuuu]\ar'[rr][rrr]\ar[uuur]
&
&
& \sr{\sone(\O,\O)}{\cov(\M)}\ar[uuuuu]
\\
  \sr{\sone(\BO,\BG)}{\p}\ar[uuuuu]\ar[rr]\ar[ur]
&
& \sr{\sone(\BO,\BO)}{\cov(\M)}\ar[uu]\ar[ur]\ar[rrr]
&
&
& \sr{\sone(\B,\B)}{\cov(\M)}\ar[uuuuu]\ar[ur]
}$$
}
\caption{}\label{survive}
\end{figure}
Almost all implications which do not appear in Figure \ref{survive} where
refuted by counter-examples (which are in fact sets of real numbers)
in \cite{coc1, coc2, CBC}.
The only unsettled implications in this diagram are marked with
dotted arrows.
We thus have the following \emph{classification problems}
(which appear in \cite{coc2} as Problems 1 and 2).

\begin{prob}
Is $\ufin(\Gamma,\Omega)=\sfin(\Gamma,\Omega)$?
\end{prob}

\begin{prob}
And if not, does $\ufin(\Gamma,\Gamma)$ imply
$\sfin(\Gamma,\Omega)$?
\end{prob}

\end{document}